\documentclass[11pt,a4paper,reqno]{amsart}
\usepackage{amsmath}
\usepackage{amsfonts}
\usepackage{amssymb}
\usepackage[utf8]{inputenc}
\usepackage[T1]{fontenc}
\usepackage{amssymb,amsmath,amsthm}
\usepackage{graphicx}
\usepackage{enumerate}
\usepackage[shortlabels]{enumitem}
\usepackage{longtable}
\usepackage{multicol}
\date{}
\newtheorem{twr}{Theorem}[section]
\newtheorem{lem}[twr]{Lemma}
\newtheorem{ques}[twr]{Question}

\title{On length of the period of the continued fraction of $n\sqrt{d}$}
\author{Filip Gawron, Tomasz Kobos}

\address{Faculty of Mathematics and Computer Science \\ Jagiellonian University \\ Lojasiewicza 6, 30-348 Krakow, Poland}

\email{filipux.gawron@student.uj.edu.pl}
\email{Tomasz.Kobos@uj.edu.pl}


\begin{document}
\maketitle
\pagestyle{empty}

\begin{abstract}
    For a given quadratic irrational $\alpha$, let us denote by $D(\alpha)$ the length of the periodic part of the continued fraction expansion of $\alpha$. We prove that for a positive integer $d$, which is not a perfect square, the sequence $(D(n\sqrt{d}))_{n=1}^{\infty}$ has infinitely many limit points.

\end{abstract}
{\centering
\section{Introduction}
}
Throughout the paper we assume that $d$ is a positive integer, which is not a perfect square. Every real number $\alpha$ has a continued fraction expansion of the form:
$$
\alpha=a_0+\cfrac{1}{a_1+\cfrac{1}{a_2+\cfrac{1}{a_3+\cdots}}}
$$
for some non-negative integers $a_1,a_2, \ldots$. We will denote this fraction by $[a_0,a_1,a_2,a_3,...]$. Lagrange showed that continued fraction of $\alpha$ is eventually periodic if and only if $\alpha$ is a quadratic irrational -- that is an irrational root of a quadratic polynomial with integer coefficients. Later Galois proved that for numbers of the form $\alpha=\sqrt{d}$ this expansion is of the form
$$\sqrt{d}=[a_0,\overline{a_1;a_2,...,a_2,a_1,2a_0}],$$
where overlined part is periodic and consists of a palindrome with added $2a_0$ at the end. By $D(\alpha)$ we denote the length of the periodic part of the continued fraction expansion of $\alpha$.

Continued fraction expansions of the quadratic irrationals and the length of their periodic parts have been widely studied by many authors. In the paper from 1972, Chowla and Chowla \cite{Ch} asked the following question.

\begin{ques}\label{ques1}
For a given integer $k \geq 1$, are there infinitely many integers $d \geq 1$ such that $D(\sqrt{d})=k$?
\end{ques}
The answers turns out to be positive, which was shown by Friesen in \cite{Fr}. In fact, Friesen's approach gives more than that. He proved that for palindromes $(a_1,a_2,\ldots, a_2, a_1)$, satisfying certain quite general condition related to the parity, there exist infinitely many integers $d \geq 1$ with $\sqrt{d}=[a_0,\overline{a_1,\ldots a_1,2a_0}]$. This was later refined by Halter-Koch, who proved that it is possible to impose some additional conditions on $d$, related to its $p$-adic valuations (see Theorem 2 in \cite{HK}). The equation $D(\sqrt{d})=k$ for small $k \geq 1$ was studied in \cite{BH} by Balkov\'a and Hru\u skov\'a. 

A similar but different line of research has been recently taken by Rada and Starosta in \cite{RS}. They studied the behaviour of the length of the continued fraction expansion of a certain transformation of $\sqrt{d}$. A \emph{M\"obius transformation} is a transformation of the form
$$h(x) = \frac{ax+b}{cx+d},$$
where $a$, $b$, $c$, $d$ are integers. Rada and Starosta were able to provide some lower and upper bounds on $D(h(\sqrt{x}))$, in terms of $D(\sqrt{x})$. These estimates are expressed with help of the determinant $|ad-bc|$. 

The main goal of the paper is to study a natural question, that falls somewhere in between the two previously mentioned lines of research. A transformation $x \to nx$, where $n \geq 1$ is an integer, is the simplest example of a M\"obius transformation. It is therefore natural to consider a variant of the Question \ref{ques1} in the class of the number of the form $n \sqrt{d}$. We shall state it in a rather general form.

\begin{ques}
\label{ques2}
For given integers $k \geq 1$ and $d \geq 1$, are there infinitely many integers $n \geq 1$ such that $D(n\sqrt{d})=k$?
\end{ques}

In the full generality, the answer to this question turns out easily to be negative, as opposed to Question \ref{ques1}. For example, if $D(\sqrt{d})$ is even, then $D(n\sqrt{d})$ is also even for every $n \geq 1$ (see Lemma \ref{even}). Therefore, in order to make this question more specific and interesting, we define 
$$A_d=\{k\in\mathbb{N}: \ \text{there exist infinitely many $n$ for which $D(n\sqrt{d})=k$}\}.$$ 
In other words, $A_d$ is the set of limit points of the sequence $(D(n\sqrt{d}))_{n=1}^{\infty}$. Our main result goes as follows.
\begin{twr}
\label{glowne}
Let $d$ be a positive integer, which is not a perfect square. Then, the set $A_d$ is infinite. In other words, the sequence $(D(n\sqrt{d}))_{n=1}^{\infty}$ has infinitely many different limit points.
\end{twr}

It does not seem possible to deduce Theorem \ref{glowne} from any of the previously mentioned results. The proof is based on establishing a lower bound on $D(n_i \sqrt{d})$ for some specific sequence $(n_i)_{i=1}^{\infty}$. In order to do this, we explore some connections between the continued fraction expansion, the Euclidean algorithm and the Pell equation. 

The paper is organized as follows. In Section \ref{lemmas} we prove some preparatory results. Theorem \ref{glowne} is proved in Section \ref{dowod}. Section \ref{problemy} concludes the paper with some further open questions.

\section{Auxiliary results}
\label{lemmas}

For each continued fraction $[a_0,a_1,a_2,a_3, \ldots]$ we define sequences $(p_i)_{i=0}^{\infty}$, $(q_i)_{i=0}^{\infty}$ as follows:\begin{multicols}{2}
\noindent
$$
\left\{ \begin{array}{ll}
p_0=a_0\\
p_1=a_0a_1+1\\
p_k=a_kp_{k-1}+p_{k-2}
\end{array} \right.
$$
$$
\left\{ \begin{array}{ll}
q_0=1\\
q_1=a_1\\
q_k=a_kq_{k-1}+q_{k-2}
\end{array} \right.
$$
\end{multicols}
Let us recall (see Chapter 1.5 and 1.6 in \cite{co}) that these sequences satisfy the equations
\begin{equation}
\label{ciag1}
    p_kq_{k-1}-p_{k-1}q_k=(-1)^{k-1}
\end{equation}
\begin{equation}\label{rat}
    \frac{p_k}{q_k}=[a_0,a_1,...,a_k]
\end{equation}
for all $k\geq 1$. The reader is referred to \cite{co} for some basics on the continued fraction expansion. 

We start our investigation with some preparatory lemmas related to the Pell equation.
Let $(x_n, y_n)_{n=1}^{\infty}$ be the sequence of the solutions of the Pell equation $x^2-dy^2=1$. It is easy to see that for every integer $a$ the sequence $x_n \pmod{a}$ is periodic mod $a$. By $m_d(a)$ we denote the period length of this sequence modulo $a$. We recall a folklore result.

\begin{lem}
\label{pell}
For any prime $p>2$ we have $m_d(p)|p^2-1$.
\end{lem}

In the next two lemmas we relate the Pell equation to $D(n \sqrt{d})$. 
\begin{lem}
\label{okres2}
For every positive integer $d$, there exists a positive integer $n$ such that $D(n \sqrt{d})=2.$
\end{lem}
\begin{proof}
We will show that, if a pair $(a,b)$ fulfills a Pell equation
$$a^2-db^2=1,$$
then 
$$b\sqrt{d}=[a-1,\overline{1,2(a-1)}].$$
Let $\beta=[\overline{1,2(a-1)}]$. We have
$$\beta=[1,2(a-1),\beta]=\frac{p_1\beta+p_0}{q_1\beta+q_0}=\frac{(2a-1)\beta+1}{2(a-1)\beta+1}$$
Solving this for $\beta$ we obtain:
$$\beta=\frac{(a-1)+\sqrt{a^2-1}}{2(a-1)},$$
as we can exclude the second negative root. Therefore
$$a-1+\frac{1}{\beta}=a-1+\frac{2(a-1)}{(a-1)+\sqrt{a^2-1}}=\frac{a^2-1+(a-1)\sqrt{a^2-1}}{(a-1)+\sqrt{a^2-1}}=\sqrt{a^2-1}=b\sqrt{d},$$
which concludes the proof.
\end{proof}
\begin{lem}
\label{even}
Suppose that $D(\sqrt{d})$ is even. Then $D(n \sqrt{d})$ is even for every $n \geq 1$.
\end{lem}
\begin{proof}
It is known (see for example \cite[Theorem 7.26]{hs}), that a negative Pell equation
$$x^2-dy^2=-1,$$
has a solution (and thus infinitely many) if and only if, the number $D(\sqrt{d})$ is odd. Therefore, if $D(n \sqrt{d})$ is odd, then the equation
$$x^2-(dn^2)y^2=-1$$
has an integer solution $(a,b)$. But then, the equation
$$w^2-dz^2=-1$$
has an integer solution $(a,nb)$. This implies that $D(\sqrt{d})$ is odd and gives us a contradiction.
\end{proof}

Now we turn our attention to the Euclidean algorithm. For given positive integers $x, y$, by $L(x, y)$ we denote the length of the Euclidean algorithm for $x$ and $y$. Here we assume that for every positive integer $n \geq 1$ we have $L(n, 1)=1$. So, for example, $L(25, 7) =4$, because the Euclidean algorithm in this case goes as follows: $(25, 7) \to (7, 4) \to (4, 3) \to (3, 1)$. We suppose that for the Fibonacci sequence $(F_n)_{n=0}^{\infty}$ we have $F_0=0$, $F_1=F_2=1$ and $\varphi$ denotes the golden ratio. A following lemma gives a lower bound on $L(a, b)$ for $a$ and $b$ with ratio close to $\varphi$.

\begin{lem}
\label{fib}
Let $a>b$ be positive integers such that
$\left | \frac{a}{b} - \varphi \right | < \left | \frac{F_{k+2}}{F_{k+1}} - \varphi \right |.$
Then $L(a, b) \geq k$.
\end{lem}
\begin{proof}

Let $x_0=a$, $x_1=b$ and for $n \geq 2$ let $x_n$ be the smaller of the numbers obtained in the $(n-1)$-th step of the Euclidean algorithm for $a$ and $b$, i.e. $x_n$ is the remainder of $x_{n-2}$ modulo $x_{n-1}$. We shall prove inductively that $x_n = (-1)^{n}aF_{n-1} + (-1)^{n+1}bF_{n}$ for $1 \leq n \leq k$. The statement is clearly true for $n=1$, so let us assume that it is true for some $1 \leq n \leq k-1$. We will show that $2x_{n} > x_{n-1}$. Indeed, by the definition of $x_n$, it is equivalent to
$$2(-1)^{n}F_{n-1}a + 2(-1)^{n+1}F_{n}b > (-1)^{n+1}F_{n-2}a + (-1)^{n}F_{n-1}b$$
or
$$(-1)^n(2F_{n-1} + F_{n-2})a > (-1)^{n}(2F_n + F_{n-1})b.$$
Simplifying we get
$$ (-1)^nF_{n+1}a > (-1)^{n} F_{n+2}b,$$
which is finally equivalent to
$$(-1)^n\frac{a}{b} > (-1)^{n}\frac{ F_{n+2}}{F_{n+1}} $$
and this clearly follows from our assumption and the monotonic convergence of $\frac{F_{n+2}}{F_{n+1}}$ to $\varphi$. In particular, $x_k > 2x_{k-1}$ and therefore $x_k \neq 0$. Hence $L(a, b) \geq k$ and conclusion follows.
\end{proof}
In the last lemma of this section we establish some relations concerning the middle of part of the palindrome with its right end.
\begin{lem}
\label{wlasnosci}
Suppose that $k=2l=D(\sqrt{d})$ is even. Then, the following properties are true:
\begin{equation}
\label{wlasnosci1}
q_{k-1} =q_{l-1}(q_l+q_{l-2})=q_{l-1}(a_l q_{l-1} + 2q_{l-2}),
\end{equation}
\begin{equation}
\label{wlasnosci2}
p_{k-1} = a_0 q_{k-1} + q_{k-2},
\end{equation}
\begin{equation}
\label{wlasnosci3}
q_{l-1}|q_{k-2}+(-1)^{l-1},
\end{equation}
\begin{equation}
\label{wlasnosci4}
a_lq_{l-1} + 2q_{l-2}|q_{k-2} + (-1)^l.
\end{equation}
\end{lem}
\begin{proof}
We will use an equivalent matrix definition of $(p_i)_{i=0}^{\infty}, (q_i)_{i=0}^{\infty}$, which follows directly from construction of these sequences (see \cite{Poo} for details). For 
$$\sqrt{d}=[a_0,\overline{a_1,...,a_{l-1},a_l,a_{l-1},...,a_1,2a_0}]$$
we have
$$
\begin{pmatrix}
p_{k-1} & p_{k-2}\\
q_{k-1} & q_{k-2}
\end{pmatrix}=
\begin{pmatrix}
a_0 & 1 \\
1 & 0 
\end{pmatrix}
\cdot
\ldots`
\cdot
\begin{pmatrix}
a_{l-1} & 1 \\
1 & 0 
\end{pmatrix}
\cdot
\begin{pmatrix}
a_l & 1 \\
1 & 0 
\end{pmatrix}
\cdot
\begin{pmatrix}
a_{l-1} & 1 \\
1 & 0 
\end{pmatrix}
\cdot
\ldots
\cdot
\begin{pmatrix}
a_1 & 1 \\
1 & 0 
\end{pmatrix}
=$$
$$=
\begin{pmatrix}
a_0 & 1 \\
1 & 0 
\end{pmatrix}
\cdot
\ldots
\cdot
\begin{pmatrix}
a_l & 1 \\
1 & 0 
\end{pmatrix}
\cdot
\left(
\begin{pmatrix}
a_0 & 1 \\
1 & 0 
\end{pmatrix}
\cdot
\ldots
\cdot
\begin{pmatrix}
a_{l-1} & 1 \\
1 & 0 
\end{pmatrix}\right)^T
\cdot
\begin{pmatrix}
a_0 & 1 \\
1 & 0 
\end{pmatrix}^{-1}
.$$
Thus
$$
\begin{pmatrix}
p_l & p_{l-1} \\
q_l & q_{l-1} 
\end{pmatrix}
\cdot
\begin{pmatrix}
p_{l-1} & p_{l-2} \\
q_{l-1} & q_{l-2} 
\end{pmatrix}^T=
\begin{pmatrix}
p_{k-1} & p_{k-2}\\
q_{k-1} & q_{k-2}
\end{pmatrix}
\cdot
\begin{pmatrix}
a_0 & 1 \\
1 & 0 
\end{pmatrix}
.$$
Finally, we obtain:
$$
\begin{pmatrix}
p_{l-1}(p_l+p_{l-2}) & p_{l}q_{l-1}+p_{l-1}q_{l-2} \\
q_{l}p_{l-1}+q_{l-1}p_{l-2} & q_{l-1}(q_l+q_{l-2}) 
\end{pmatrix}=
\begin{pmatrix}
a_0p_{k-1}+p_{k-2} & p_{k-1} \\
a_0q_{k-1}+q_{k-2} & q_{k-1} 
\end{pmatrix}
.$$
By comparing matrix entries we get:
\begin{itemize}
\item Property $(\ref{wlasnosci1}) $, by using the equality $q_l=a_lq_{l-1}+q_{l-2}.$.
\item Property $(\ref{wlasnosci2})$: by comparing matrix entries we have $p_{k-1}=p_lq_{l-1}+p_{l-1}q_{l-2}$. Then, using $(\ref{ciag1})$ we can write it as
$$q_lp_{l-1}+(-1)^{l-1}+q_{l-1}p_{l-2}+(-1)^{l-2}=q_lp_{l-1}+q_{l-1}p_{l-2}=a_0q_{k-1}+q_{k-2}.$$
\item Property $(\ref{wlasnosci3})$: from previous properties we have
$$p_{k-1}=p_lq_{l-1}+p_{l-1}q_{l-2}=a_0q_{k-1}+q_{k-2}.$$
From $(\ref{wlasnosci1})$ we get:
$$a_0q_{k-1}+q_{k-2}=a_0q_{l-1}(a_lq_{l-1}+2q_{l-2})+q_{k-2}.$$
Comparing this two equalities we obtain:
$$a_0q_{l-1}(a_lq_{l-1}+2q_{l-2})+q_{k-2}=p_lq_{l-1}+p_{l-1}q_{l-2}=q_{l-1}p_{l}+q_{l-1}p_{l-2}+(-1)^{l-2}.$$
Therefore
$$q_{k-2}+(-1)^{l-1}=q_{l-1}(p_{l}+p_{l-2}-a_0(a_lq_{l-1}+2q_{l-2})).$$
\item Property $(\ref{wlasnosci4})$ similarly to property $(\ref{wlasnosci3})$:
$$a_0q_{l-1}(a_lq_{l-1}+2q_{l-2})+q_{k-2}=p_lq_{l-1}+p_{l-1}q_{l-2}=q_{l}p_{l-1}+(-1)^{l-1}+p_{l-1}q_{l-2}$$
$$q_{k-2}+(-1)^{l}=(a_lq_{l-1}+2q_{l-2})(p_{l-1}-a_0q_{l-1}),$$
where we used the relation $q_l+q_{l-2}=a_lq_{l-1}+2q_{l-2}$.

\end{itemize}

\end{proof}
It is worth mentioning that property (\ref{wlasnosci2}) was used also by Friesen \cite{Fr} in the proof of the affirmative answer to Question \ref{ques1}. 
\section{Proof of Theorem \ref{glowne}}
\label{dowod}

Before proving our main result, we need one more lemma, which connects the number $D(n\sqrt{d})$, a corresponding Pell equation and the length of the Euclidean algorithm (as promised in the Introduction). In fact, Theorem \ref{glowne} is an easy consequence of this lemma.

\begin{lem}
\label{nieparzyste}
Let $p, q$ be odd numbers. If $2|D(\sqrt{d})$ and $(x_a,y_a)$ is a solution of a Pell equation $x^2-dy^2=1$,  satisfying the conditions:
\begin{itemize}
    \item $2pq|y_a$,
    \item $p|x_a\pm1$,
    \item $q|x_a\mp1$,
\end{itemize}
then there exist infinitely many positive integers $n$ such that 
$$D(n\sqrt{d}) \in \{ 2L(p,q), 2(L(p,q)+1) \}.$$
\end{lem}
\begin{proof}
If there exists at least one solution satisfying above conditions, then there exist infinitely many of them, as the sequence of solutions of Pell's equation is periodic modulo $2pq$. Let us take one of these solutions $(x_m,y_m)$, such that $x_m, y_m > 4p^2q^2$. We can consider another Pell equation of the form
$$w^2 - d\left ( \frac{y_m}{2pq} \right )^2z^2 = 1$$
Clearly $(w, z) = (x_m, 2pq)$ is a solution of this equation. In fact, it is the fundamental solution. Indeed, let $(w_1, z_1)$ be the fundamental solution and $(w_2,z_2)$ be the next one. Assume that $z_1 < 2pq$. Then 
$$2pq \geq z_2 = 2w_1z_1 \geq 2w_1 > \frac{y_m}{2pq},$$
which contradicts our choice of $m$. Therefore $(w_1, z_1)=(x_m, 2pq)$. Let $n = \frac{y_m}{2pq}$ and $k = D(n \sqrt{d})$. By Lemma \ref{even} we know that $k$ is even, so let us write $k=2l$. From \cite{co} (see Chapter 4.8) we know that the fundamental solution $(w_1,z_1)$ is a pair $(p_{k-1},q_{k-1})$. Hence
\begin{equation}
\label{pellpara}
    (x_m,2pq)=(p_{k-1},q_{k-1})
\end{equation}
First part of Lemma \ref{wlasnosci} yields the equality 
\begin{equation}
\label{rownanie}
2pq = q_{l-1}(a_l q_{l-1} + 2q_{l-2}).
\end{equation}
Also from Lemma \ref{wlasnosci} we easily get the following congurences
$$p_{k-1}\equiv (-1)^{l} \pmod{q_{l-1}},$$
$$p_{k-1}\equiv (-1)^{l-1} \pmod{a_lq_{l-1} + 2q_{l-2}}.$$
Using $(\ref{pellpara})$ we can write this in the form:
$$q_{l-1}|x_m+(-1)^{l},$$
$$a_lq_{l-1} + 2q_{l-2}|x_m+(-1)^{l-1}.$$
Since $p|x_m\pm1$, we have that $\gcd(p,a_lq_{l-1} + 2q_{l-2})=1$ or $\gcd(p,q_{l-1})=1$ (as $p$ is odd). Similar property holds also for $q$. As $p$ and $q$ divide $x_m+(-1)^\alpha$ with different parity of $\alpha$, the two inequalities: $\gcd(p,a_lq_{l-1} + 2q_{l-2})>1$ and $\gcd(q,a_lq_{l-1} + 2q_{l-2})>1$ can not be true at the same time. Without loss of generality, we assume that 
$$\gcd(p,a_lq_{l-1} + 2q_{l-2})=1.$$ Then, from $(\ref{rownanie})$ it follows that $p|q_{l-1}$. Furthermore, if $2|q_{l-1}$, then $2|a_lq_{l-1} + 2q_{l-2}$, so again by $(\ref{rownanie})$ we have $4|2pq$, which is false. In consequence, $2|a_lq_{l-1} + 2q_{l-2}$, which implies that $2|a_l$ and 
$$p=q_{l-1},$$
$$q=\frac{a_l}{2}q_{l-1} + q_{l-2}.$$

Now, we consider the number $L(p,q)$. We have that 
$$L(p,q)=L\left ( q_{l-1},\frac{a_l}{2}q_{l-1}+q_{l-2} \right )=L(q_{l-1},q_{l-2})+1=L(a_{l-1}q_{l-2}+q_{l-3},q_{l-2})+1=L(q_{l-2},q_{l-3})+2$$
$$=L(a_{l-2}q_{l-3}+q_{l-4},q_{l-3})+2=L(q_{l-3}, q_{l-4}) + 3= \ldots = L(q_{2},q_{1}) + l-2.$$
It follows that for $q_1=1$ we have $L(p, q) = l-1$ and for $q_1>1$ we have $L(p, q)=l$.
So we proved that for $n=\frac{y_m}{2pq}$ we have 
$$D(n\sqrt{d}) \in \{ 2L(p,q), 2(L(p,q)+1) \}.$$
Since we have infinitely many options for choosing sufficiently large $m$, this finishes the proof.
\end{proof}

Now we are ready to prove our main result.

\emph{Proof of Theorem \ref{glowne}}. 
By Lemma \ref{okres2} it is enough to consider the case $D(\sqrt{d})=2$ -- every limit point of the sequence $D(nc\sqrt{d})$ for a fixed $c$, is also a limit point of the sequence $D(n\sqrt{d})$.

Let $r$ be any odd positive integer. Suppose that there exists an odd prime $p$ not dividing $d$, but dividing $x_{8r}+1$ where $(x_{8r},y_{8r})$ is the $8r$-th solution of the Pell equation
$$x^2-dy^2=1.$$
Let $t$ be the positive integer such that $p \in (\phi^{2t},\phi^{2(t+1)}).$

The distance between two consecutive numbers of the form $\frac{a}{p}$ is equal to $\frac{1}{p}$. Thus, we can choose a positive integer $b$ such that
$$\left|\frac{b}{p}-\phi\right|<\frac{1}{p}<\phi^{-2t}.$$
We note also that
$$ \left|\frac{F_{t+1}}{F_t}-\phi\right|=\frac{\phi+\phi^{-1}}{\phi^{2t}+(-1)^t}>\phi^{-2t}.$$
Therefore, by Lemma \ref{fib} we have
$$L(b,p)\geq t-1.$$
Now, let $q$ be any prime such that:
\begin{enumerate}[a)]
\item $2r$ and $\frac{q^2-1}{8}$ are relatively prime
\item $q \equiv b \pmod{p}$.
\end{enumerate}
Since $8|q^2-1$ the system of congruences 
$$m \equiv 8r \pmod{16r}, \quad m \equiv 0 \pmod{q^2-1}$$
has a solution $m$.
Let us recall that $p|x_{8r}+1$ and $p$ does not divide $d$. Hence, $p$ divides $y_{8r}$. Thus, looking modulo $p$ the pair $(x_{8r}, y_{8r})$ is congruent to the pair $(-1, 0)$ and therefore the pair $(x_{16r}, y_{16r})$ is congruent to $(1, 0)$. This shows that $m_d(p)|16r$. Hence, from the first congruence it follows that $p|x_m+1$. From the second congruence and Lemma \ref{wlasnosci} we get that $m_d(q)|m$ and thus $q|x_m-1$. This shows that $p$ and $q$ satisfy conditions of Lemma \ref{nieparzyste}. It follows that for infinitely many positive integers $n$ we have:
$$D(n\sqrt{d}) \in \{ 2L(p,q), 2(L(p,q)+1) \}.$$
From $q \equiv b \pmod{p}$ we obtain
$$L(p,q)=L(b,p)\geq t-1=\left \lfloor \frac{\log_{\phi}p}{2} \right \rfloor - 1.$$
On the other hand, we have
$$L(p,q)\leq \log_{\phi}p+1.$$
This shows that for infinitely many $n$ we have
$$D(n\sqrt{d}) \in (\log_{\phi}p-3,2\log_{\phi}p+4).$$
In particular, we get a limit point in each interval $(\log_{\phi}p-3,2\log_{\phi}p+4)$.

To finish the proof, we are left with proving that we can choose infinitely many primes $p$ for a variable $r$. As $d$ has only finitely many divisors, it is enough to prove that there are infinitely many primes $p$ such that $p|x_{8r}+1$ for some $r \geq 1$. Let us consider $r$ prime. In this case, we have $m_d(p)|16r$. If $r$ does not divide $m_d(p)$, then $m_d(p)|16$ which gives us $p|y_{16}$ and therefore it is satisfied by only finitely many primes $p$. If $r|m_d(p)$, then by Lemma \ref{pell} we have also that $r|p^2-1$. This clearly yields the desired infinite number of possibilities for $p$ and the proof is finished. \qed

\section{Concluding remarks}
\label{problemy}

It is natural to ask, if something more can be said in general about the set $A_d$ of the limit points of the sequence $(D(n\sqrt{d}))_{n=1}^{\infty}$, besides the fact that it is of infinite cardinality. More specifically, we pose the following question.

\begin{ques}\label{hip1}
Is it true, that for every non-square integer $d \geq 1$ and every $k \geq 1$ at least one of the numbers $k, k+1$ belongs to $A_d$?
\end{ques}

Such a conjecture may seem to be quite strong at first glance, but there is a motivation behind it. We recall that by $L(x, y)$ we denoted the length of the Euclidean algorithm applied for $x$ and $y$. Let us state perhaps somewhat more natural question, that is related directly to $L(x, y)$.

\begin{ques}\label{hip2}
Is it true, that for every integer $k \geq 1$, there exists an integer $N \geq 1$, such that for every $n > N$ and $1 \leq i \leq k$ there exists $1 \leq m \leq n$ such that $L(m, n)=i$?
\end{ques}

It turns out, that a slight modification of the argument used in the proof of Theorem \ref{glowne} shows that an affirmative answer to Question \ref{hip2} would directly imply an affirmative answer to Question \ref{hip1}. We do not know if any of these questions has a positive answer, but we believe that some further properties of $A_d$ could be established in the full generality.

\end{document}